\documentclass[12pt]{amsart}
\usepackage{color}
\usepackage{epsfig}
\usepackage{verbatim}
\usepackage{parskip}

\usepackage{algorithm}
\usepackage{algpseudocode}
\usepackage{amssymb}

\newtheorem{theorem}{Theorem}[section]
\newtheorem{lemma}[theorem]{Lemma}

\begin{document}

\title {On cyclically 4-connected cubic graphs}
\maketitle

\small
 
\begin{center}
R. J. Kingan \footnote{Bloomberg LP., 
120 Park Avenue,
New York, NY 10165,
rkingan@bloomberg.net.}, 
S. R. Kingan \footnote{Department of Mathematics,
Brooklyn College, City University of New York,
 Brooklyn, NY 11210,
skingan@brooklyn.cuny.edu.} 
 
\end{center}

\begin{abstract} 
For $k \ge 4$, let $Q_{2k}$ and $V_{2k}$ denote the ladder and M{\"o}bius ladder on $2k$ vertices, respectively. We prove results that build on a result by Wormald that states that any cyclically $4$-connected cubic graph other than $Q_8$ or $V_8$ is obtained from a smaller cyclically $4$-connected cubic graph by bridging a pair of non-adjacent edges. We introduce the concept of cycle spread, which generalizes the edge pair distance defined by Wormald, and show that the set of pairs of edges that needs to be considered in order to obtain all cyclically $4$-connected cubic graphs is smaller than the set of all pairs of non-adjacent edges. We prove that all non-planar cyclically $4$-connected cubic graphs with at least $10$ vertices, other than the M{\"o}bius ladders and the Petersen graph, are obtained from $Q_8$ by bridging pairs of edges with cycle spread at least $(1,2)$. Moreover every graph obtained in this way is non-planar, cyclically $4$-connected, and cubic. All planar cyclically $4$-connected cubic graphs with at least $10$ vertices except for the ladders are obtained from the ladders by bridging pairs of edges with cycle spread at least $(1,2)$. We implemented an algorithm based on these results using McKay's nauty system for isomorphism checking.

\end{abstract}

\section{Introduction}\label{Introduction}

A cubic graph is a graph whose vertices have degree 3. If $n$ is the number of vertices and $m$ is the number of edges, then $n$ is even and $m=\frac{3n}{2}$. A $k$-connected graph is one in which at least $k$ vertices must be removed to disconnect the graph. 

The degree of every vertex in a cubic graph is 3, so it cannot have connectivity higher than 3. As a result a concept called cyclic $k$-connectivity is used to define finer levels of connectivity. An {\it edge-cut} is a set of edges whose removal disconnects the graph. An edge-cut $X$ is {\it cycle-separating} if at least two of the components in the disconnected graph formed when $X$ is removed have cycles. A 3-connected cubic graph is {\it cyclically $k$-connected} if it has  at least $n\ge 2k$ vertices and it has no cycle-separating edge-cut of size less than $k$. Moreover, a cyclically $k$-connected cubic graph will not have any cycles of length less than $k$ since if it does, the $k$ edges incident to the $k$ vertices of the cycle, but not in the cycle, may be removed to obtain two components with cycles. In this paper, we present structural results for cyclically $4$-connected cubic graphs.

A cycle with $k$ vertices is called a  {\it $k$-cycle}. A cycle with 3 vertices is called a {\it triangle} and a graph without triangles is called {\it triangle-free}. A cyclically 4-connected cubic graph is triangle-free, has $n\ge 8$ vertices, and when removal of a set of 3 edges results in a disconnected graph, only one component has cycles. 

\begin{figure}[h]
\centering
\includegraphics[width=1.3in]{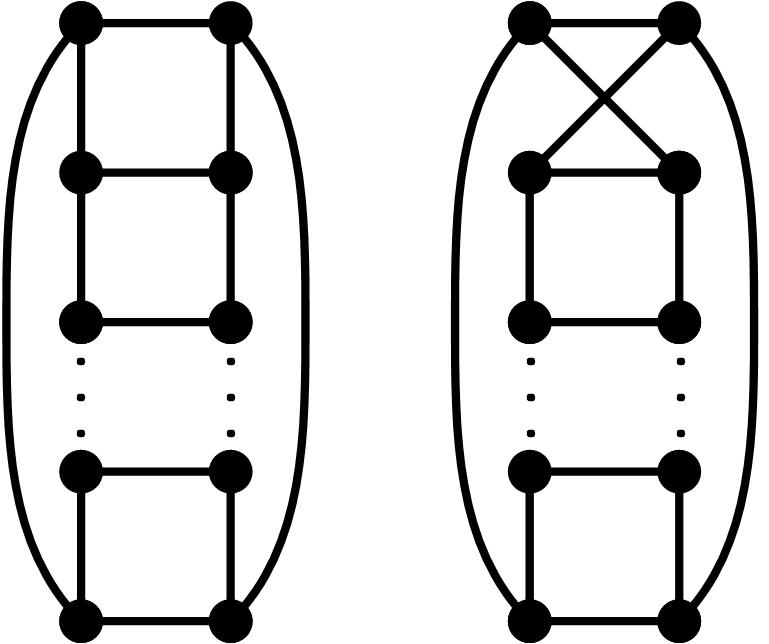}
\caption{ $Q_{2k}$, $k \ge 4$ (left) and $V_{2k}$, $k \ge 3$ (right). \label{LaddersMobiusLadders}}
\end{figure}

Before proceeding, familiarity with some cyclically 4-connected cubic graphs is needed. Figure \ref{LaddersMobiusLadders} displays two well-known infinite families: the ladders $Q_{2k}$ for $k\ge 4$ and the M{\"o}bius ladders $V_{2k}$ for $k\ge 4$. It is easy to see that $Q_8$ and $V_8$ are the only 8-vertex cyclically 4-connected cubic graphs. Further, note that $V_6 \cong K_{3,3}$ and it is not cyclically $4$-connected since it has fewer than 8 vertices. 

\begin{figure}[h]
\centering
\includegraphics[width=3.5in]{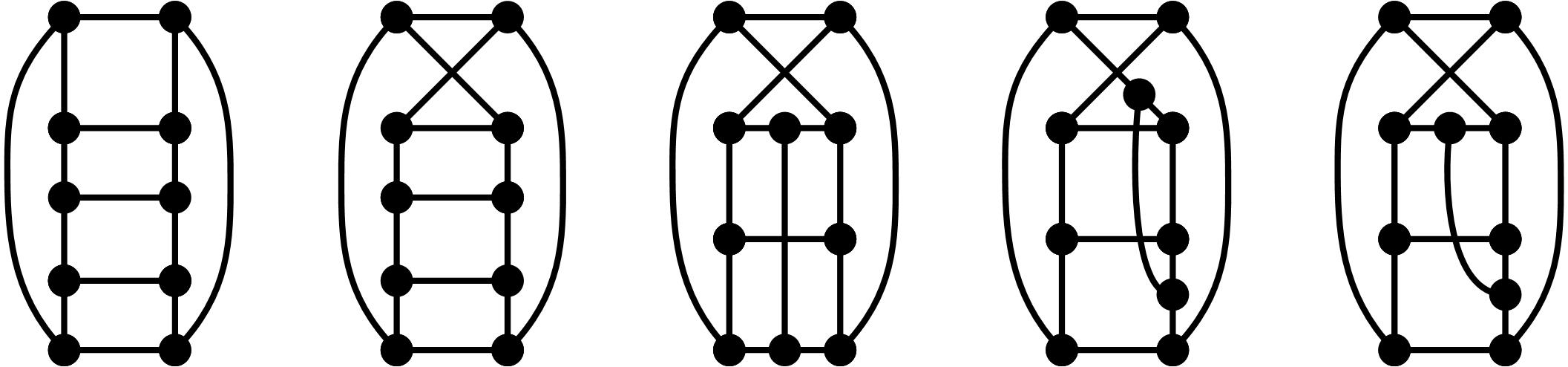}
\caption{ The non-isomorphic cyclically 4-connected cubic graphs with 10 vertices. \label{10-vertex-C4C-Cubicgraphs}}
\end{figure}

The five graphs shown in Figure \ref{10-vertex-C4C-Cubicgraphs}, the first three of which are $Q_{10}$, $V_{10}$, and the Petersen graph $P_{10}$, are the only cyclically 4-connected cubic graphs on 10 vertices. The Online Encyclopedia of Integer Sequences at {\tt http://www.oeis.org} gives the numbers of 3-connected (A204198) and cyclically 4-connected (A175847) cubic graphs for $n\le 22$ vertices.  

\begin{figure}[h]
\centering
\includegraphics[width=2in]{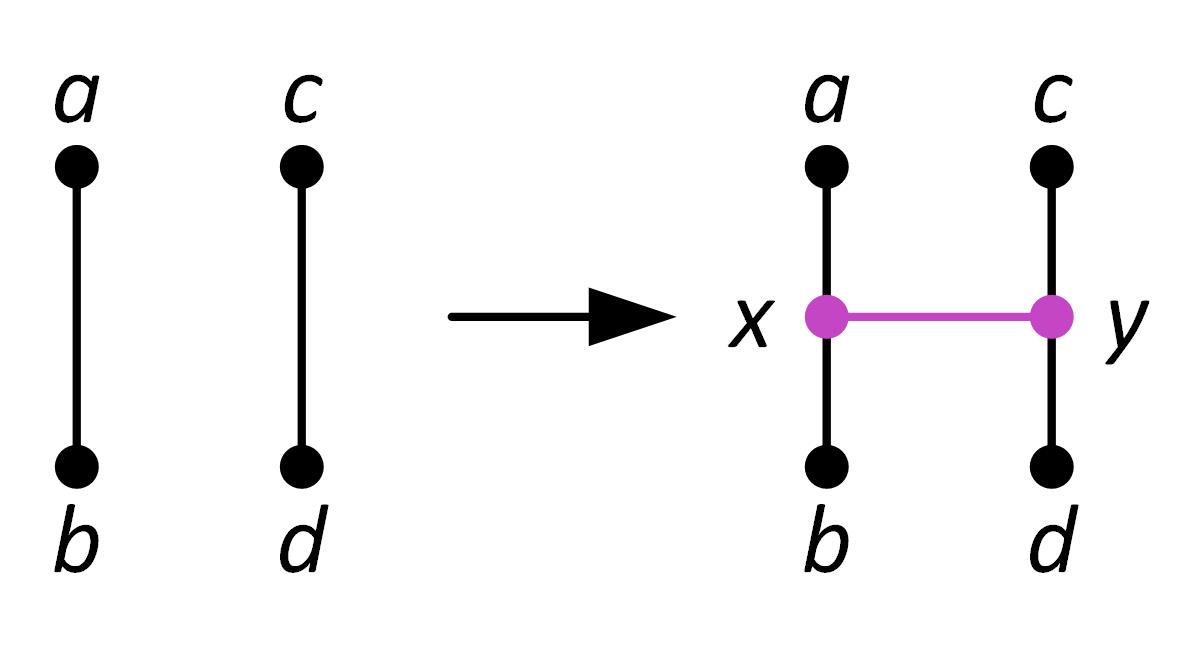}
\caption{ Bridging two edges.  \label{BG-2}}
\end{figure}

Two distinct vertices are {\it adjacent} if they are joined by an edge, and {\it non-adjacent} otherwise. Two distinct edges are {\it adjacent} if they have a common vertex, and {\it non-adjacent} otherwise. A pair of distinct edges $ab$ and $cd$ in a graph $G$ is \textit{bridged} if they are subdivided by vertices $x$ and $y$, respectively, forming paths of length 2, and $x$ and $y$ are joined by an edge as shown in Figure \ref{BG-2}.  The new graph $G'$, called a {\it bridge addition} of $G$, has two new vertices $x$ and $y$ and edges $xa$, $xb$, $yc$, $yd$, and $xy$. Observe that $G'$ has three more edges than $G$. We say $G'$ is obtained from $G$ by bridging a pair of edges.

Let $G$ be a graph and $e=uv$ be an edge with end vertices $u$ and $v$. The graph with edge $e$ deleted is called an {\it edge-deletion} and is denoted by $G\backslash e$ or~$G\backslash uv$. When deleting edge $e$, the end vertices $u$ and $v$ remain. To contract edge $e$, collapse the edge by identifying the end vertices $u$ and $v$ as one vertex, and delete the resulting loop. The graph with edge $e$ contracted is called an {\it edge-contraction} and denoted by  $G / e$.  
A graph $H$ is   a {\it minor} of a~graph $G$ if $H$ can be obtained from $G$ by deleting edges (and any isolated vertices formed as a result) and contracting~edges. We write $H=G\backslash X/Y$, where $X$ is the set of edges deleted and $Y$ is the set of edges contracted. 

To reverse the bridging operation, we must delete the new edge $xy$ and contract edges $xb$ and $yc$ in the graph $G'$. The order in which the deletion and two contractions are done does not matter. Thus reversing bridging, which we refer to as {\it unbridging}, is equivalent to the operation $G' \backslash xy / \{xb, yc\}$.  It is well known that edge-deletion and edge-contraction preserve planarity. Since the unbridging operation can be expressed in terms of edge-deletions and edge-contractions, it also preserves planarity.
 
\begin{figure}[h]
\centering
\includegraphics[width=4.25in]{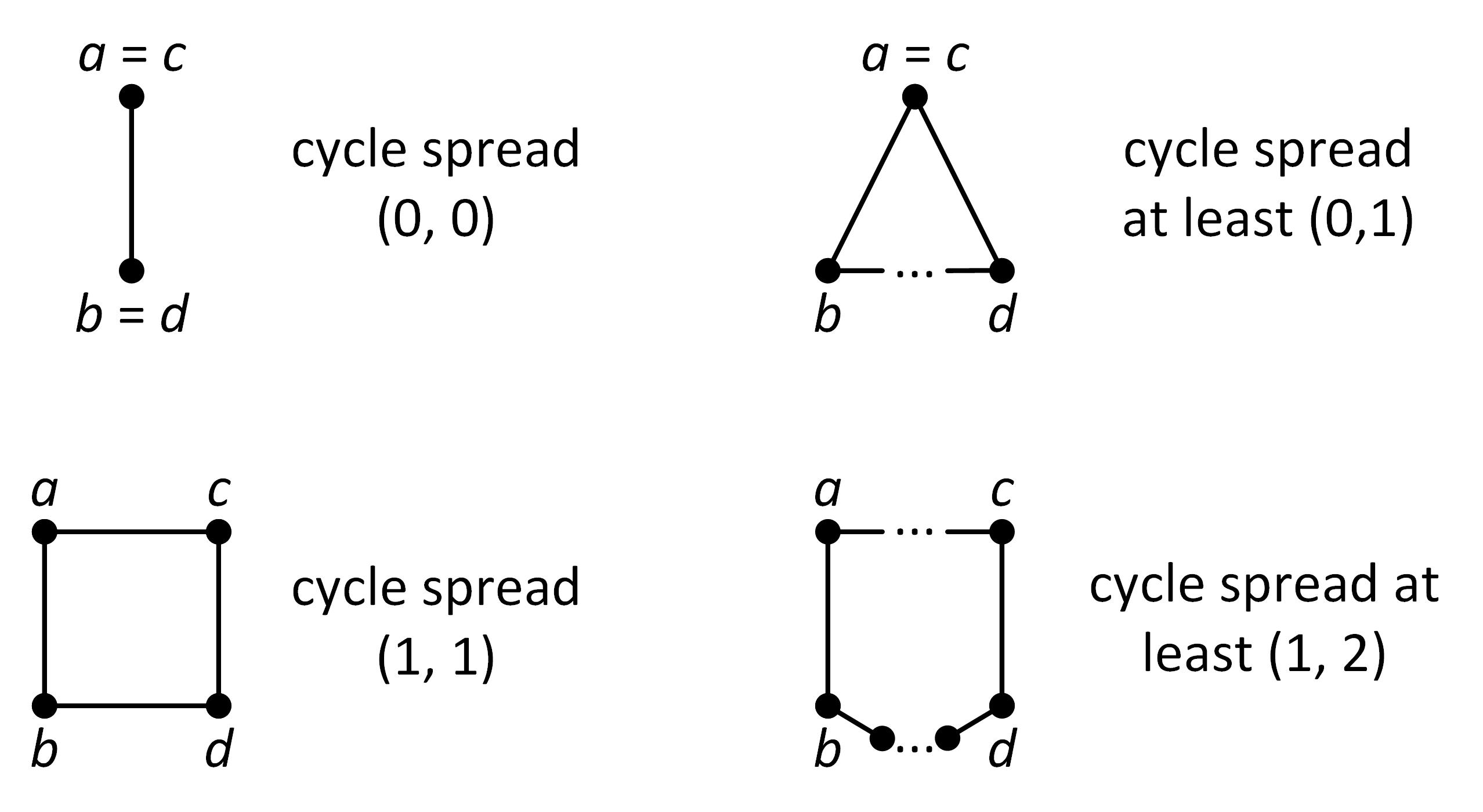}
\caption{ Cycle spread examples.  \label{cycle-spread-figure}}
\end{figure}

In a $2$-connected graph, every pair of edges is contained in at least one cycle. Let $ab$ and $cd$ be a pair of edges in a $2$-connected graph, and let $C$ be a smallest cycle containing $ab$ and $cd$. There are two distinct paths in $C$ connecting $ab$ and $cd$; without loss of generality we may assume that the two paths connect $a$ to $c$ and $b$ to $d$, respectively, and that the $a-c$ path is no longer than the $b-d$ path. We can then associate a pair $(i, j)$ where $i \le j$ with $ab$ and $cd$ that consists of the lengths of these two paths. We call this pair the {\it cycle spread} of edges $ab$ and $cd$.   Figure \ref{cycle-spread-figure} shows four examples of cycle spread: 

\begin{enumerate}
\item By convention, if edges $ab$ and $cd$ are the same, their cycle spread is $(0, 0)$.

\item If $ab$ and $cd$ are adjacent but not identical, then their cycle spread is $(0, j)$ for some $j > 0$.

\item If $ab$ and $cd$ are non-adjacent but part of a $4$-cycle, then their cycle spread is $(1, 1)$

\item If $ab$ and $cd$ are non-adjacent and not part of any $4$-cycle, then their cycle spread is $(i, j)$ for some $i > 0, j > 1$. In terms of lexicographic ordering on pairs of integers, we can say that their cycle spread is at least $(1, 2)$.

\end{enumerate}  

In this paper, we will prove that if $G$ and $G'$ are triangle-free $2$-connected cubic graphs such that $G$ can be obtained from $G'$ by bridging two edges with cycle spread at least $(1, 1)$, then with the exception of two families of graphs, namely $Q_{2k}$ and $V_{2k}$, $G$ can be obtained from $G'$ by bridging two edges with cycle spread at least $(1, 2)$. As a consequence, we will prove that:

\begin{itemize}

\item Any non-planar cyclically $4$-connected cubic graph is obtained from $Q_8$ by repeatedly bridging pairs of edges with cycle spread at least $(1, 2)$; and 

\item Any planar cyclically $4$-connected cubic graph is obtained from $Q_{2k}$ for some $k \ge 4$ by repeatedly bridging pairs of edges with cycle spread at least $(1, 2)$.

\end{itemize}

A standard approach to graph generation is to generate a large collection of graphs, and then discard the graphs that are not in the desired set. Our goal is to generate precisely the desired set.

Section 2 contains preliminaries and the significance of our results. Sections 3 and 4 contain the proofs of the main theorems. In Section 5, we will discuss how these results can be used to develop an algorithm to exhaustively and correctly generate cyclically $4$-connected cubic graphs, and how a simple modification of the approach can also be used to generate cyclically $5$-connected cubic graphs, although the latter procedure is not exhaustive.


\section {Preliminaries}\label{Preliminaries}

The wheel graph with $n$ vertices and $n-1$ spokes is denoted by $W_{n-1}$, where $n\ge 4$.  Let $G$ and $H$ be 3-connected cubic graphs such that $H$ is a minor of $G$ and $G\not\cong W_3$.  Tutte proved that $G$ can be obtained from $H$ by repeatedly bridging edges. At each stage the graph obtained remains 3-connected and cubic \cite{Tutte1967}. 

Thus every $3$-connected cubic graph $G$, other than $W_3$, with $n$ vertices and $m$ edges can be obtained from a smaller $3$-connected cubic graph $G'$ with $m - 3$ edges and $n - 2$ vertices by bridging a pair of edges. This process begins with $W_3$. The non-isomorphic graphs obtained by bridging a pair of distinct edges in $W_3$ are $K_{3,3}$ and the prism graph (the geometric dual of $K_5 \backslash e$). The only non-isomorphic graph obtained by bridging edges in $K_{3,3}$ is $V_8$. The only non-isomorphic graph obtained by bridging edges in the prism graph is $Q_8$.

Wormald proved that every cyclically $4$-connected cubic graph except $Q_8$ and $V_8$ can be obtained from a cyclically $4$-connected cubic graph with fewer vertices by bridging non-adjacent edges \cite{Wormald1979}; that is, pairs of edges with cycle spread at least $(1, 1)$. Note that he referred to cyclically 4-connected graphs as ``irreducible graphs.'' 

\begin{theorem} \label{Wormald} Let $G$ be a cyclically $4$-connected cubic graph such that $G\not \cong Q_8, V_8$. Then $G$ is obtained from a cyclically $4$-connected cubic graph with fewer vertices by bridging pairs of edges with cycle spread at least $(1,1)$. Moreover, each graph obtained in this manner is a cyclically $4$-connected cubic graph.
\end{theorem} 

Note that Wormald did not specify that cyclically $4$-connected graphs must have $n \ge 8$ vertices. He considered $V_6 \cong K_{3,3}$ to be cyclically $4$-connected. In the statement of his theorem (Theorem 4 in \cite{Wormald1979}), he requires $G$ to be a cyclically $4$-connected graph such that $G \not \cong W_3, Q_8$. Our definition of cyclically $4$-connected is from Robertson, Seymour, and Thomas \cite{RST-C5C}. In particular, we do not consider $K_{3,3}$ to be cyclically $4$-connected and we start the process from $8$ vertices.

Wormald defines the {\it distance} between two edges $e$ and $f$ in a connected graph to be one less than the length of the shortest path which includes $e$ and $f$. By this definition, two adjacent edges have distance 1, corresponding to cycle spread $(0, j)$ for some $j \ge 1$, and two non-adjacent edges have distance at least 2, corresponding to cycle spread at least $(1, 1)$. Further, he suggested that bridging two non-adjacent edges of distance at least $k-2$ in a cyclically $k$-connected graph produces another cyclically $k$-connected graph. He added  that this appeared to be  a promising approach to the investigation of cyclically $k$-connected graphs, but that the analogy breaks down for $k\ge 5$. Graphs obtained by bridging two non-adjacent edges in a cyclically 5-connected  graph of distance at least 3 will be cyclically 5-connected, but this procedure does not give all the cyclically $5$-connected graphs. He gave an infinite family of cyclically 5-connected graphs that cannot be obtained in this manner. 

The {\it chromatic number} of $G$ is the smallest number of colors required so that adjacent vertices in the graph are assigned different colors. The {\it edge chromatic number} of $G$ is the smallest number of colors required so that adjacent edges in the graph are assigned different colors. Vizing proved that the edge-chromatic number of a graph with maximum degree $\Delta$ is either $\Delta$ or $\Delta + 1$ \cite{Vizing1964}. Thus cubic graphs have edge-chromatic number either 3 or 4.  A {\it snark} is a bridgeless cubic graph with edge-chromatic number 4.  A {\it non-trivial snark} is cyclically 4-connected, to avoid simple examples. Snarks were named by Martin Gardener because, like the snarks in Lewis Carrol's poem ``Hunting the snarks,'' they are elusive creatures. The Petersen graph is the smallest snark. It is well known that snarks are non-planar. One of the most significant unsolved problems in this area is to construct precisely snarks from smaller snarks, if such an algorithm is even possible. Our results are a step in this direction because they lead to an algorithm for generating all non-planar cyclically $4$-connected cubic graphs.

It is well known that if $G$ is a cyclically $4$-connected cubic graph with no $V_8$ minor, then $G$ is planar (assuming $K_{3,3}$ is not considered cyclically $4$-connected). A graph $H$ is {\it topologically contained} in a graph $G$, if a subdivision of $H$ is isomorphic to a subgraph of $G$. Kelmans \cite{Kelmans1993} and independently Mahary \cite{Mahary2000} proved that if $G$ is a cyclically 4-connected cubic  graph that does not topologically contain $Q_8$, then $G$ is isomorphic to $P_{10}$ or $V_{2k}$, where $k\ge 4$.  This result can be obtained from our main results by relating deletion and contraction to the bridging operation. In fact, due to the connection we made between bridging and minors, in future work we will extend our results to minimally $3$-connected cyclically $4$-connected graphs. In \cite{Kingan2021}, we already developed theorems and an algorithm for constructing minimally $3$-connected graphs. We will use the concept of cycle spread to eliminate graphs that are not cyclically $4$-connected.

\section{Triangle-free $2$-connected cubic graphs}

The main theorem of this section shows that for a triangle-free $2$-connected graph, under certain conditions, a bridge added to a pair of edges with cycle spread $(1,1)$ can be ``exchanged'' for a bridge added to a pair of edges with cycle spread at least $(1,2)$.

\begin{theorem} (Exchange Theorem) \label{exchange-theorem} Let $H$ be a triangle-free $2$-connected cubic graph with $n \ge 8$ vertices, and let $G$ be obtained from $H$ by bridging a pair of edges with cycle spread at least $(1,1)$. Then one of the following statements is true:

\begin{enumerate}

\item[i)] $H \cong Q_{2k}$ and $G \cong Q_{2(k + 1)}$ for some $k \ge 4$;

\item[ii)] $H \cong V_{2k}$ and $G \cong V_{2(k + 1)}$ for some $k \ge 4$; or

\item[iii)] $G$ is obtained from $H$ by bridging a pair of edges with cycle spread at least $(1,2)$.

\end{enumerate}

\end{theorem}

\begin{proof} Since $H$ is a triangle-free $2$-connected cubic graph, and since $G$ is obtained from $H$ by bridging a pair of edges with cycle spread at least $(1,1)$, $G$ will also be a triangle-free $2$-connected cubic graph. Let $ab$ and $cd$ denote the pair of edges with cycle spread at least $(1,1)$ in $H$ that can be bridged to obtain $G$. 

Next, since $H$ is $2$-connected, $ab$ and $cd$ are in a cycle. However, since $H$ has no triangles, $ab$ and $cd$ are not part of any triangle. If $ab$ and $cd$ have cycle spread at least $(1,2)$, then there is nothing to prove. Therefore suppose that $ab$ and $cd$ have cycle spread exactly $(1,1)$. Then $ab$ and $cd$ form a 4-cycle with two other edges; without loss of generality say that $ac$ and $bd$ are also edges. Thus we have the 4-cycle $acdba$, as shown in Figure \ref{bridge-lemma-2}.


\begin{figure}[h]
\centering
\includegraphics[width=1.65in]{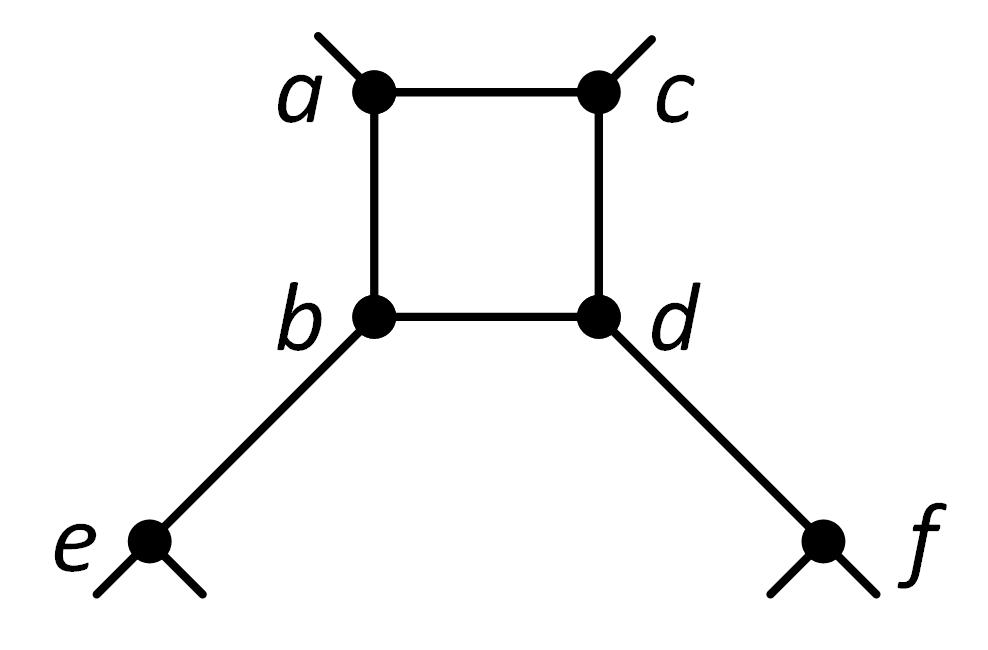}
\caption{Bridging a pair of edges with cycle spread $(1,1)$ with edges $be$ and $df$ incident to $b$ and $d$ \label{bridge-lemma-2}}
\end{figure}

Now $H$ is cubic, therefore $b$ and $d$ are each incident with exactly one other edge outside the 4-cycle. Denote these edges by $be$ and $df$ as shown in Figure \ref{bridge-lemma-2}. Note that edges $be$ and $df$ are not  adjacent, since $H$ is triangle-free. Moreover, $b$ and $d$ are incident to no other edges since each vertex has degree 3.  

Observe that the graph formed by bridging edges $ab$ and $cd$ with edge $xy$ is isomorphic to the graph formed by bridging edges $be$ and $df$ with edge $xy$, as indicated in Figure \ref{bridge-lemma-3}.

\begin{figure}[h]
\centering
\includegraphics[width=3.25in]{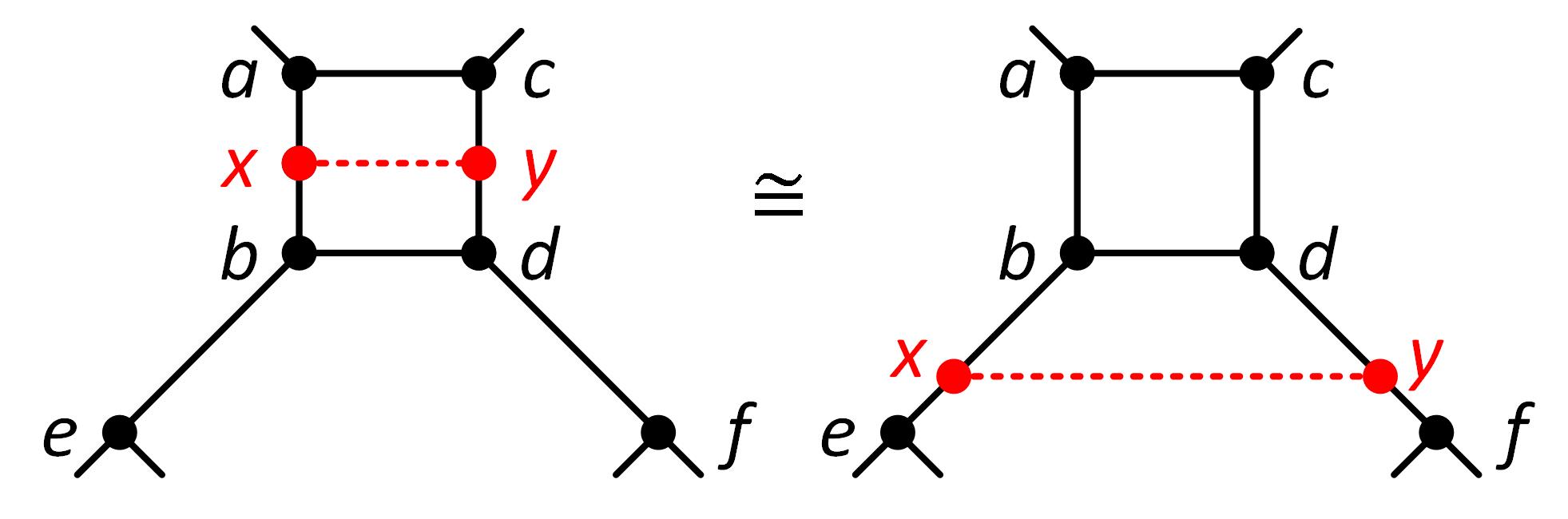}
\caption{Exchanging one pair of bridged edges for another \label{bridge-lemma-3}}
\end{figure}

Note that Figure \ref{bridge-lemma-3} illustrates just a portion of a large cubic graph. Instead of bridging edges $ab$ and $cd$, we may bridge edges $ac$ and $bd$, as shown in Figure \ref{bridge-lemma-4}, which is just Figure \ref{bridge-lemma-3} viewed sideways and relabeled. Here the bridge could be exchanged for a bridge of the edges incident with vertices $c$ and $d$ that are not part of the cycle.

\begin{figure}[h]
\centering
\includegraphics[width=2.5in]{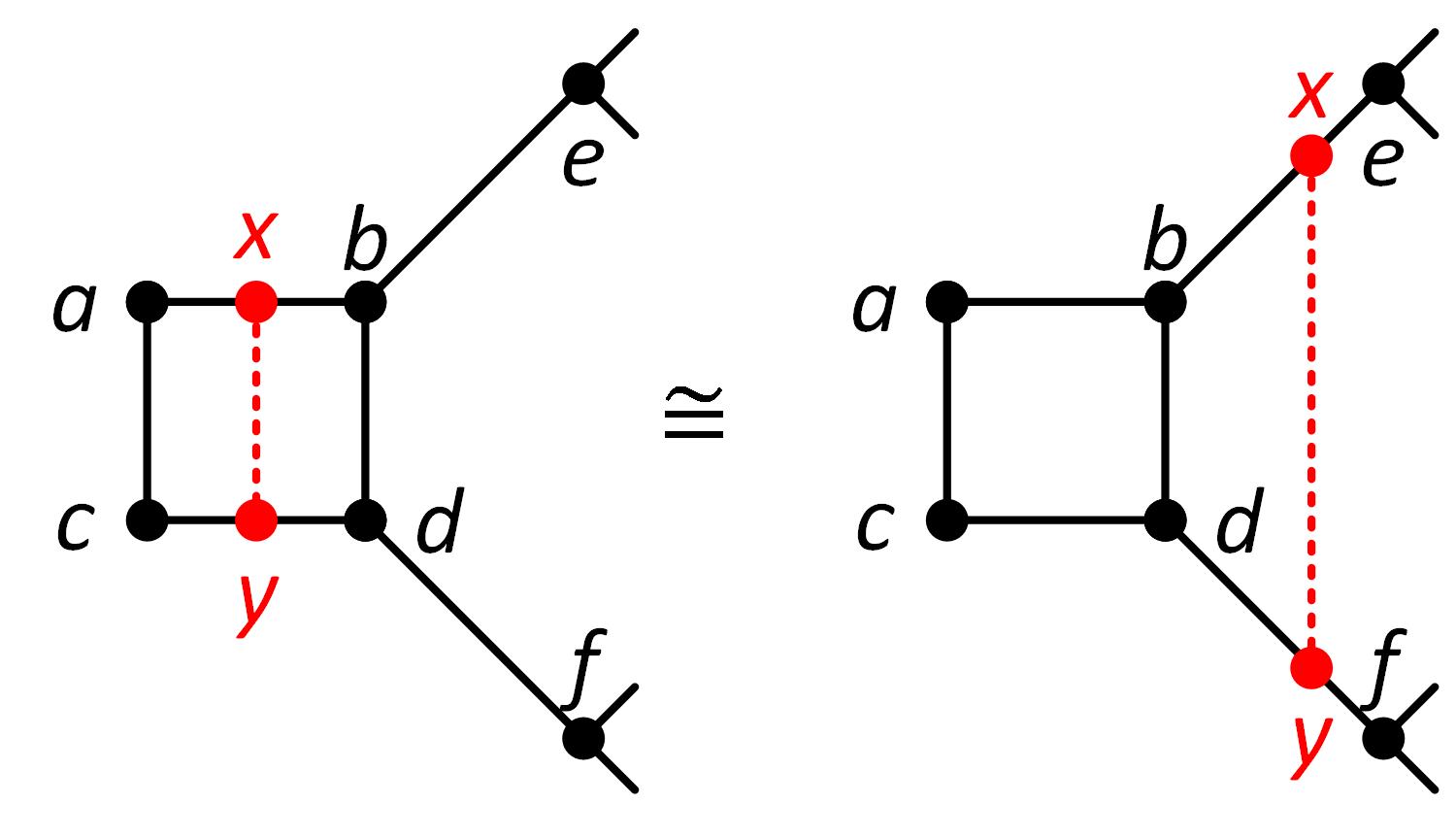}
\caption{Exchanging one pair of bridged edges for another, sideways \label{bridge-lemma-4}}
\end{figure}

If vertices $e$ and $f$ are not adjacent, then edges $be$ and $df$ have cycle spread at least $(1,2)$, and the result follows. Therefore suppose vertices $e$ and $f$ are adjacent, so that the cycle spread of edges $be$ and $df$ is exactly $(1,1)$. Then vertices $e$ and $f$ are each incident to exactly one other edge, forming the 4-cycle $bdfeb$, as shown in Figure \ref{bridge-lemma-5}. Those two edges are not adjacent because $H$ has no triangles, and we can proceed as above.  We think of the labeled edge $xy$ as hopping along the four cycles.

\begin{figure}[h]
\centering
\includegraphics[width=0.875in]{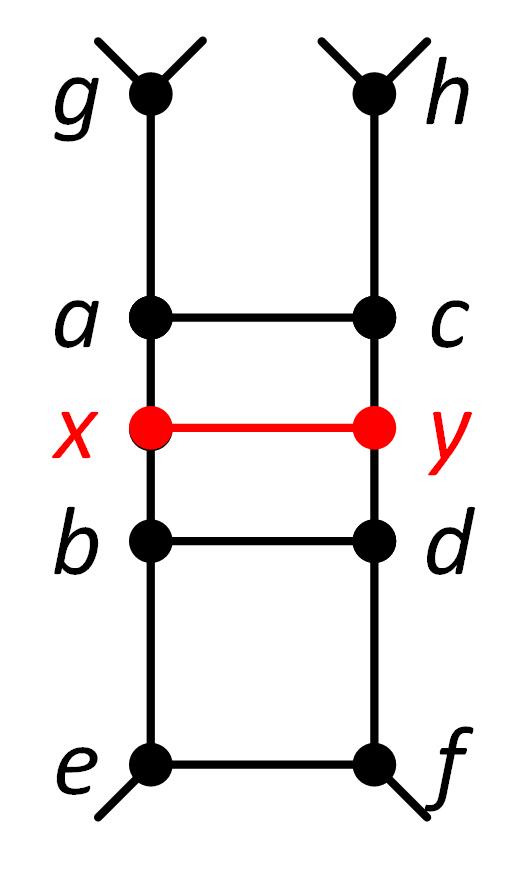}
\caption{Exchanging a bridge of one edge pair with cycle spread $(1,1)$ for another. \label{bridge-lemma-5}}
\end{figure}

As we continue hopping along 4-cycles, eventually either we will reach a stage where we are bridging two edges with cycle spread at least $(1,2)$ or we will reach the same 4-cycle $acdba$ where we began. The only way for the latter to occur is if $H \cong Q_{2k}$ or $H \cong V_{2k}$, as shown in Figure \ref{cube-mobius-exchange-figure}. When $H\cong Q_{2k}$, the edge labeled $xy$ hops along the rungs of the ladder returning to the initial 4-cycle $acdba$ in exactly the same manner. When $H\cong V_{2k}$, due to the twist that occurs edge $xy$ returns to the initial 4-cycle $acdba$ with vertex $y$ along edge $ab$ and vertex $x$ along edge $cd$. But this graph is isomorphic to the original graph; moreover with another sequence of hops we reach the original graph exactly.

\begin{figure}[h]
\centering
\includegraphics[width=5.15in]{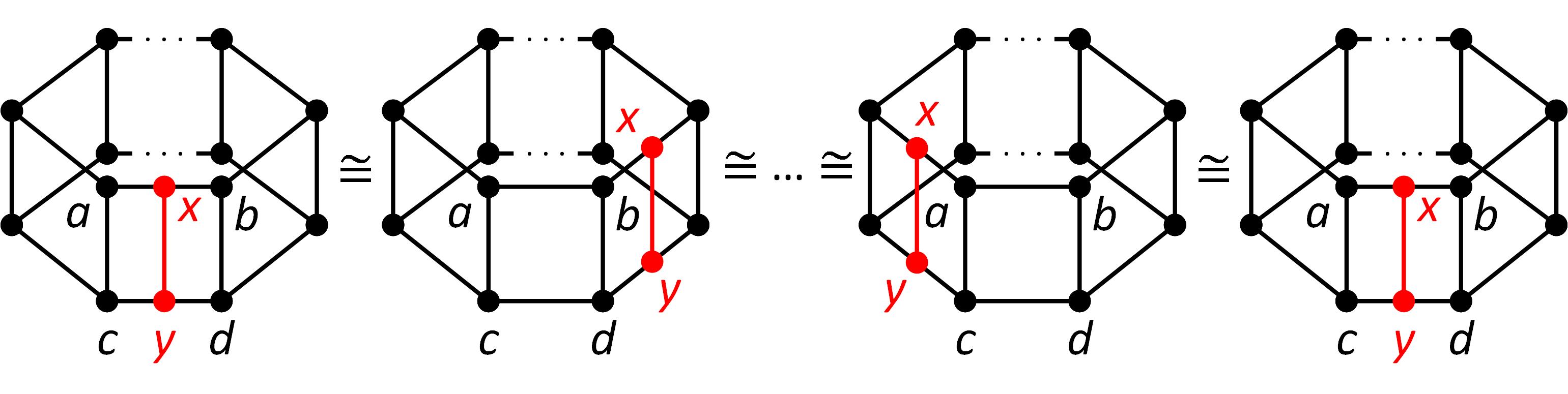}
\includegraphics[width=5.15in]{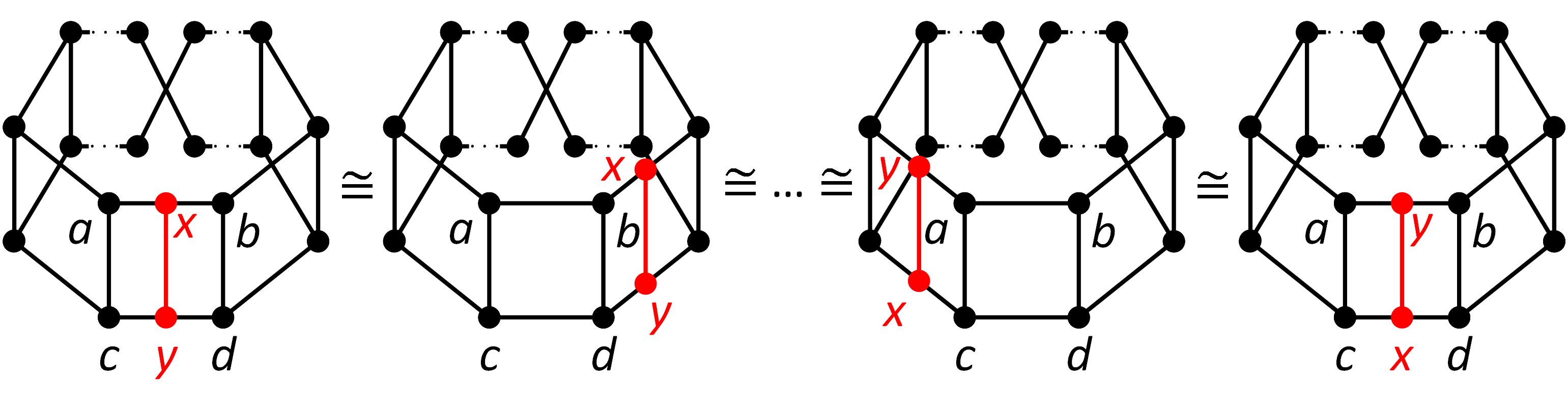}
\caption{Exchanging bridges of pairs of edges with cycle spread $(1,1)$ in $Q_{2k}$ and $V_{2k}$. \label{cube-mobius-exchange-figure}}
\end{figure}

It is apparent that these are the only ways to avoid reaching a pair of edges with cycle spread at least $(1,2)$, because each exchange of a bridge of edges $ab$ and $cd$ for a bridge of edges $be$ and $df$ requires that there be an edge $bd$. Since the graph is cubic there can be no ``forks in the road''; the sequence of exchanges proceeds along a ladder of $4$-cycles until it ends at a pair of edges with cycle spread at least $(1,2)$, or circles back on itself.
\end{proof}

The theorem is stated for all $k \ge 8$ instead of all $k \ge 6$ because the graph analogous to the ladder with six vertices is not triangle free. The only triangle-free $2$-connected cubic graph on six vertices is $K_{3,3}$, which is isomorphic to $V_6$. Every non-adjacent pair of edges in $K_{3,3}$ has cycle spread exactly $(1,1)$ and bridging any pair of non-adjacent edges in $K_{3,3}$ gives $V_8$.

It is worth noting that $H$ is only required to be a triangle-free $2$-connected graph. If $G$ is obtained from $H$ by bridging a 4-cycle, then we can traverse a ladder of 4-cycles until we either reach a pair of edges with cycle spread at least $(1,2)$, in which case we are done, or reach the 4-cycle where we started, in which case $H$ must be $Q_{2k}$ or $V_{2k}$ as indicated in Figure \ref{bridge-lemma-3}. Observe that in Figures  \ref{bridge-lemma-2} and \ref{bridge-lemma-3}, the subgraph shown is connected to the remainder of $G$ only via the additional edges incident with vertices $a$, $c$, $e$, and $f$, because the graph is cubic. This is what allows us to exchange a bridge of $ab$ and $cd$ with a bridge of $be$ and $df$ and get isomorphic graphs.

\begin{figure}[h]
\centering
\includegraphics[width=2.3in]{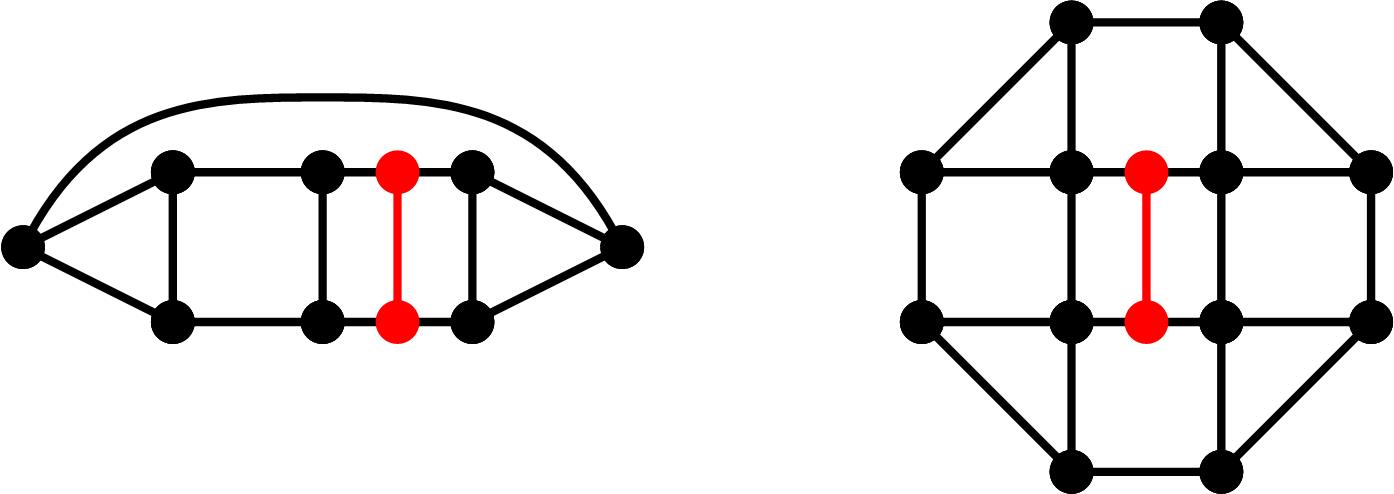}
\caption{Counterexamples to Theorem \ref{exchange-theorem} if $G$ contains triangles (left) or is not cubic (right). \label{exchange-counterexamples}}
\end{figure}

The result is not true if the graph contains triangles, as in the graph on the left in Figure \ref{exchange-counterexamples}, because the edges incident to two vertices in a 4-cycle may be adjacent. The result is also not true if the graph is not cubic as in the graph on the right in Figure \ref{exchange-counterexamples}, because then the vertices in the 4-cycle may be incident to other edges. In this case the graph produced by the exchange is not isomorphic to the graph with the original bridge. Finally, the result is not true if the graph is not $2$-connected since there may be pairs of edges that do not lie in any cycle.

\section {Structure of cyclically 4-connected cubic graphs}

In this section we will prove that we obtain all cyclically 4-connected cubic graphs, with a few exceptions, by starting with $Q_8$ and bridging only pairs of edges with cycle spread at least $(1, 2)$. Specifically, we will show that for $n \ge 10$ vertices, all  non-planar cyclically $4$-connected cubic graphs, except for the Petersen graph and the M\"obius ladders, can be generated from $Q_8$ by bridging pairs of edges with cycle spread at least $(1,2)$. Moreover, every graph obtained in this way is non-planar, cyclically $4$-connected and cubic. For $n \ge 10$ vertices, all planar cyclically $4$-connected graphs, except for the ladders, can be generated from the ladders by bridging pairs of edges with cycle spread at least $(1,2)$. Every graph generated in this way is cyclically $4$-connected and cubic, but not necessarily planar.

In other words, from the top-down perspective, we can ``trace'' any non-planar cyclically $4$-connected cubic graph other than $V_{2k}$ back to $Q_8$ via a sequence of bridges of pairs of edges with cycle spread at least $(1,2)$. A similar analysis holds for planar cyclically $4$-connected cubic graphs. However, unlike non-planar graphs, we can only ``trace'' a sequence of bridges of edges with cycle spread at least $(1, 2)$ back to $Q_{2k}$ for some $k \ge 4$. Then via a sequence of ``rung'' deletions, we can step down the ladder to $Q_8$. From the bottom-up perspective, we can grow precisely the non-planar cyclically $4$-connected cubic graphs beginning with $Q_8$. This is a step toward the overarching goal of generating precisely all snarks.

We begin with a lemma that is a consequence of Theorem \ref{exchange-theorem}.

\begin{lemma}\label{exchange-corollary} Let $H$ be a cyclically $4$-connected cubic graph with $n \ge 8$ vertices, and let $G$ be obtained from $H$ by bridging a pair of edges with cycle spread at least $(1,1)$. Then one of the following statements is true:

\begin{itemize}

\item[i)] $H \cong Q_{2k}$ and $G \cong Q_{2(k + 1)}$ for some $k \ge 4$;

\item[ii)] $H \cong V_{2k}$ and $G \cong V_{2(k + 1)}$ for some $k \ge 4$; or

\item[iii)] $G$ is obtained from $H$ by bridging a pair of edges with cycle spread at least $(1, 2)$.

\end{itemize}

\end{lemma}

\begin{proof} Cyclically 4-connected cubic graphs must have at least 8 vertices, and are $2$-connected and triangle-free. Hence the result follows from Theorem \ref{exchange-theorem}. \end{proof}

We can illustrate the exchange that occurs in Lemma \ref{exchange-corollary} by examining bridges of pairs of edges with cycle spread exactly $(1,1)$ in $Q_{10}$. There are two ways such a bridge can occur:

\begin{itemize}

\item By adding a ``rung'' to the graph; or

\item By bridging a pair of two ``rungs'' that appear together in a 4-cycle. 

\end{itemize}

\begin{figure}[h]
\centering
\includegraphics[width=2.40in]{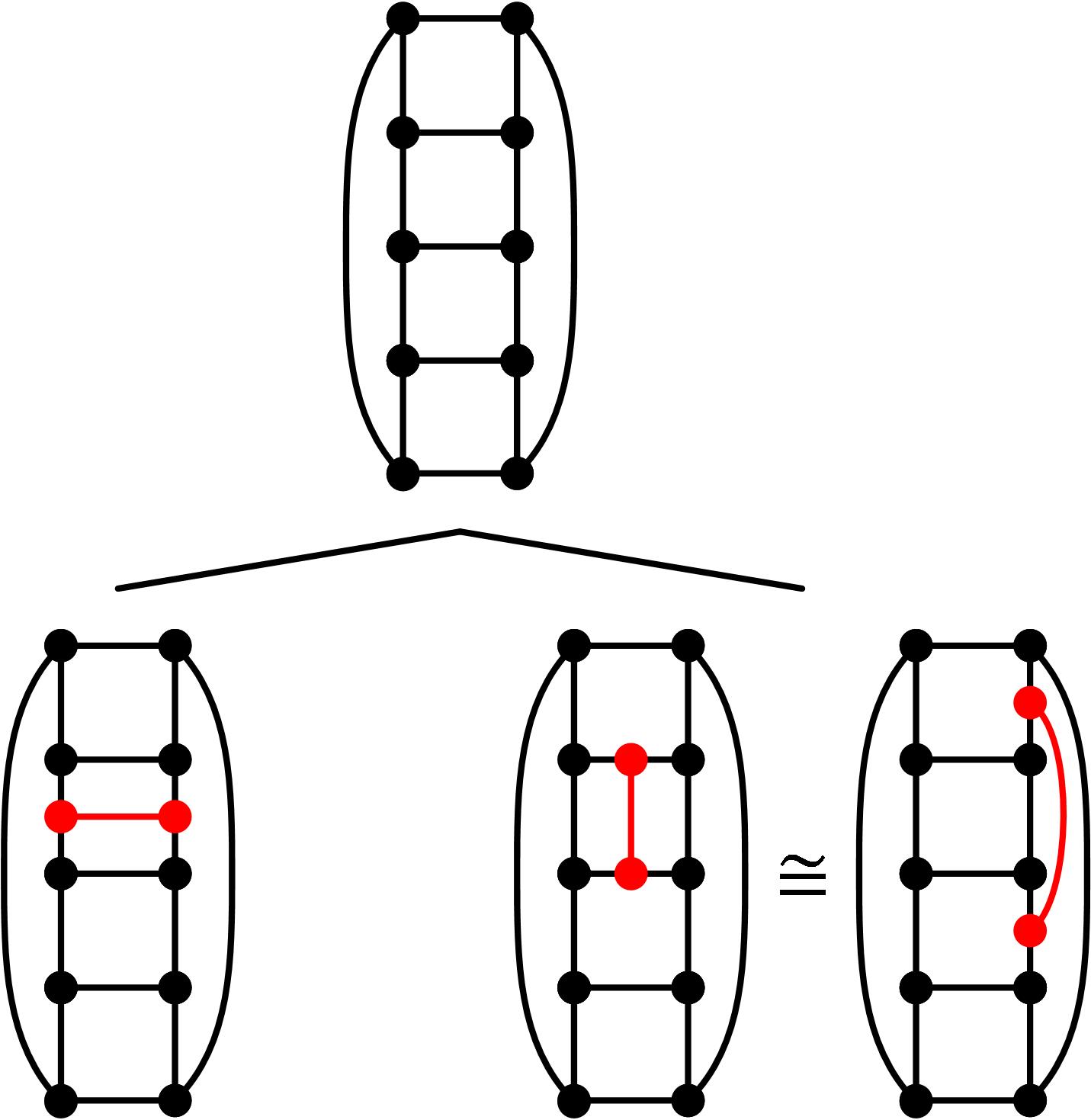}
\caption{Planar graphs obtained by bridging edges in $Q_{10}$. \label{exchange-figure}}
\end{figure}

Both ways are illustrated in Figure \ref{exchange-figure}. In the first case, shown on the left, $Q_{12}$ is produced from $Q_{10}$. In the second case, shown on the right, the bridge can be exchanged for a bridge of two edges with cycle spread $(1,2)$. 

The next result is one of the main results in this paper.

\begin{theorem} \label{MainTheorem-2} Let $G$ be a non-planar cyclically $4$-connected cubic graph such that $G\not \cong P_{10}$ or $V_{2k}$, where $k\ge 4$. Then $G$ is obtained from $Q_8$ by bridging pairs of edges with cycle spread at least $(1,2)$. Moreover, every graph obtained in this way is non-planar, cyclically $4$-connected, and cubic.
\end{theorem}

\begin{proof} The proof is by induction on $k$. Observe that $Q_8$ and $V_8$ are the only cyclically $4$-connected cubic graphs on $8$ vertices. The non-isomorphic graphs obtained by bridging pairs of edges in $Q_8$ with cycle spread at least $(1,1)$ are shown in Figure \ref{v8-q8-connection-figure}. Beginning with $Q_8$ and only bridging pairs of edges with cycle spread at least $(1,2)$, there are two distinct resulting graphs, both of which are non-planar. Bridging any pair of edges with cycle spread $(1,1)$ in $Q_8$ results in $Q_{10}$. Note that the graphs $A_{10,1}$ and $A_{10,2}$, shown in Figure \ref{v8-q8-connection-figure}, can also be obtained from $V_8$ by adding bridges to pairs of edges with cycle spread at least $(1,2)$; that is $A_{10, 2} \cong B_{10, 1}$ and $A_{10, 1} \cong B_{10, 2}$. The Petersen graph, $P_{10}$, can also be obtained from $V_8$, but not $Q_8$, this way, and $V_{10}$ can only be obtained from $V_8$ by bridging a pair of edges with cycle spread $(1,1)$. Therefore every 10-vertex graph that can be obtained from $V_8$ by bridging pairs of edges with cycle spread at least $(1,1)$, except for $P_{10}$ and $V_{10}$, can also be obtained from $Q_8$ by bridging pairs of edges with cycle spread at least $(1,2)$.

\begin{figure}[h]
\centering
\includegraphics[width=5.0in]{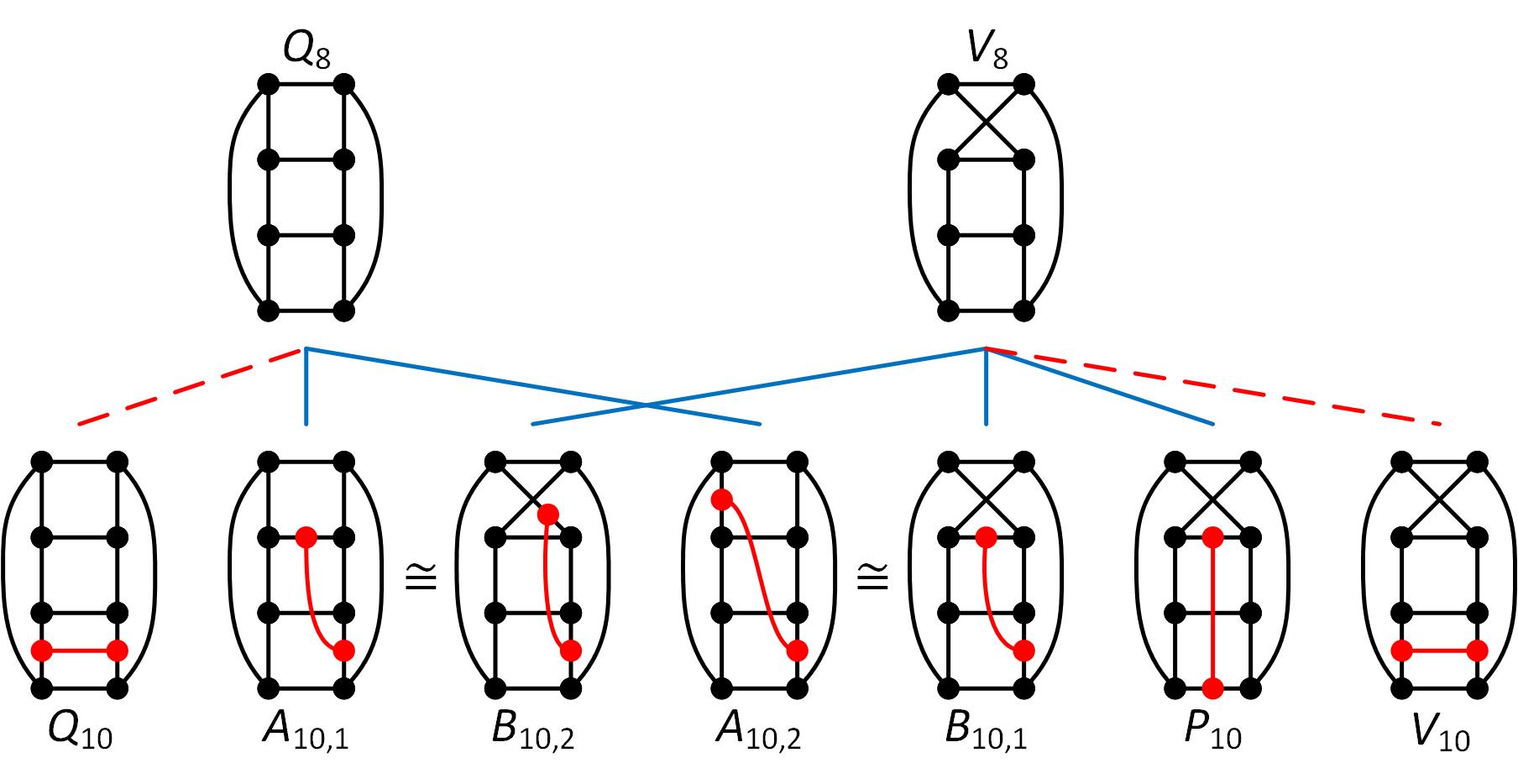}
\caption{Non-isomorphic graphs obtained by bridging pairs edges of $Q_8$ and $V_8$ with cycle spread at least $(1,1)$. Dashed lines indicate bridges that can only be added to pairs of edges with cycle spread $(1, 1)$ to produce the larger graph.\label{v8-q8-connection-figure}}
\end{figure}

In order to show that the result holds for $k = 6$, we must show that there are no non-planar cyclically $4$-connected graphs with $12$ vertices other than $V_{12}$ that cannot be obtained from $Q_8$ by bridging pairs of edges with cycle spread at least $(1,2)$. In other words we must show that any $12$-vertex graph that can be obtained from $V_{10}$ or $P_{10}$ by bridging a pair of vertices with cycle spread at least $(1,1)$ can also be obtained by bridging a pair of edges in $B_{10,1}$ or $B_{10,2}$ with cycle spread at least $(1,2)$.

Figure \ref{v10-bridges} shows the six non-isomorphic graphs other than $V_{12}$ that are obtained from $V_{10}$ by bridging pairs of edges with cycle spread at least $(1,1)$. It is easily checked that each of them is also obtained from $B_{10, 1}$ or $B_{10,2}$ by bridging pairs of edges with cycle spread at least $(1,2)$.

Likewise, Figure \ref{petersen-bridges-figure} shows the two non-isomorphic graphs that are obtained from $P_{10}$ by bridging pairs of edges with cycle spread at least $(1,1)$. It is easily checked that these two graphs are also obtained from $B_{10, 1}$ and $B_{10, 2}$ by bridging pairs of edges with cycle spread at least $(1,2)$. Therefore all 12 vertex graphs that can be obtained from $P_{10}$ by bridging pairs of edges with cycle spread at least $(1,1)$ can also be obtained from $Q_8$ by bridging pairs of edges with cycle spread at least $(1,2)$. Thus the result holds for $k=4, 5, 6$.

\begin{figure}[h]
\centering
\includegraphics[width=4.25in]{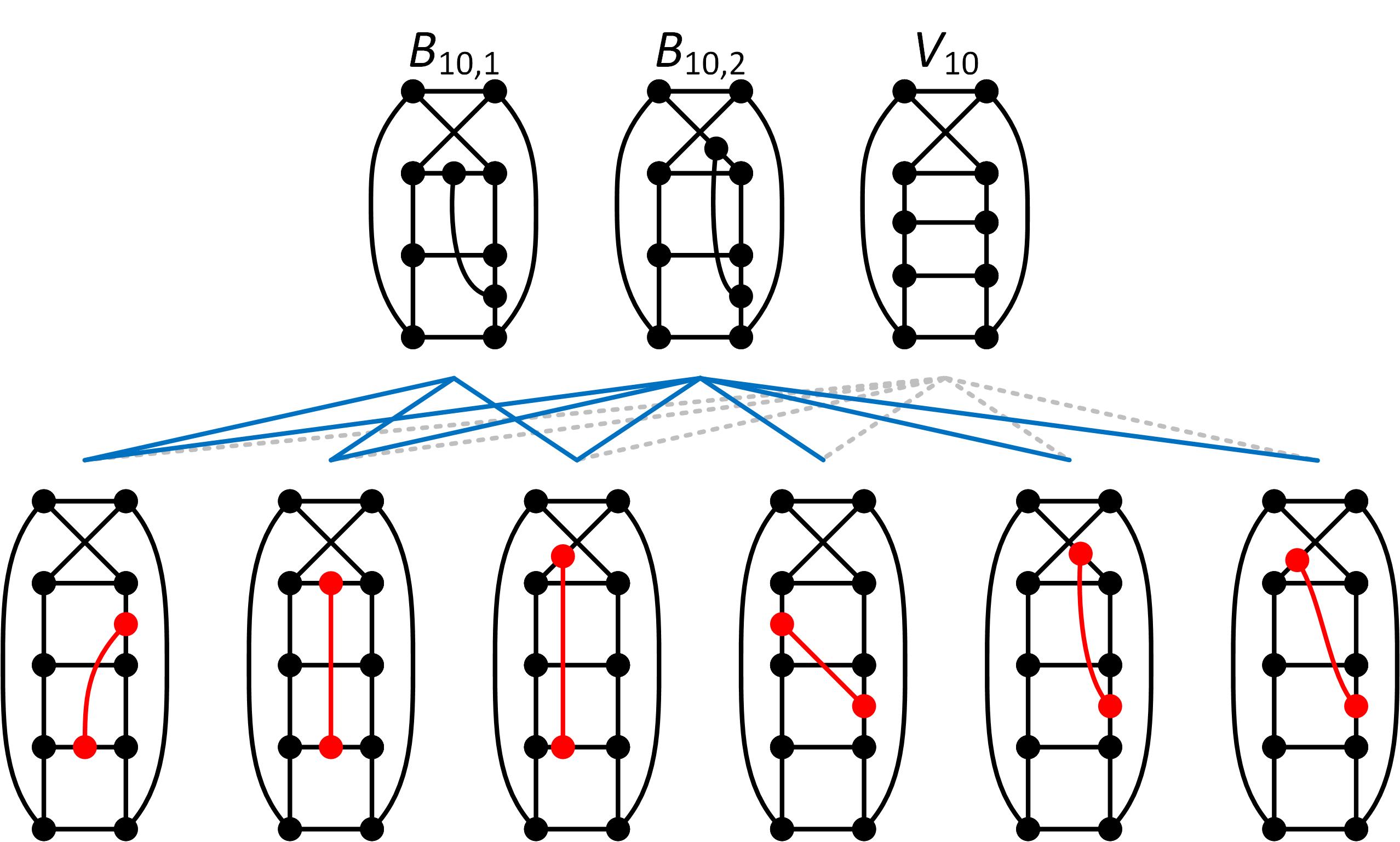}
\caption{Non-isomorphic graphs obtained by bridging pairs of edges with cycle spread at least $(1,1)$ in $V_{10}$. The first three graphs are obtained from $B_{10,1}$ by bridging a pair of edges with cycle spread at least $(1,2)$, and all of the graphs are obtained in the same way from $B_{10,2}$. \label{v10-bridges}}
\end{figure}

\begin{figure}[h]
\centering
\includegraphics[width=2.15in]{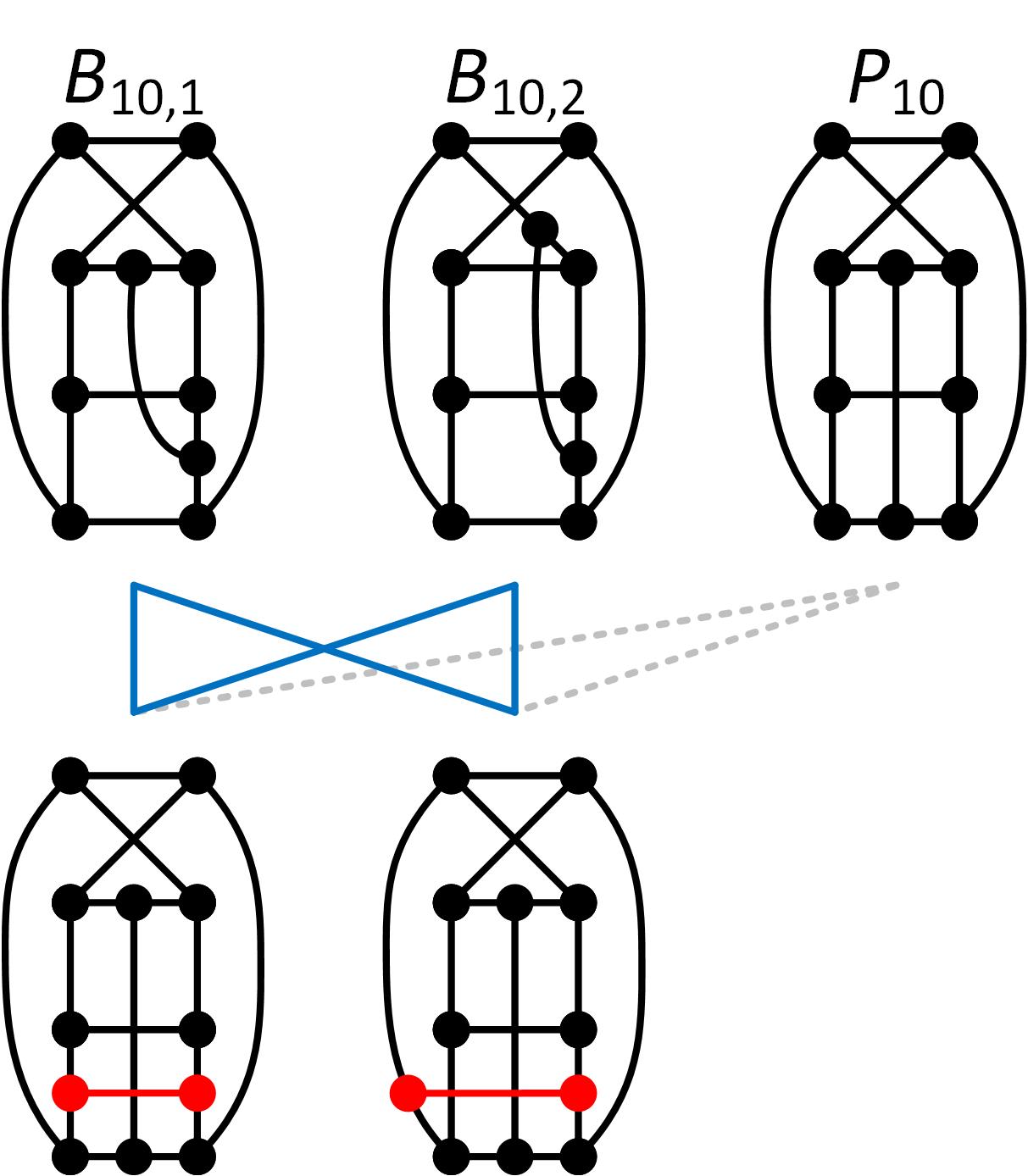}
\caption{Non-isomorphic graphs obtained by bridging pairs of edges with cycle spread at least $(1,1)$ in the Petersen graph. Both are also obtained from $B_{10,1}$ and $B_{10,2}$.\label{petersen-bridges-figure}}
\end{figure}

Suppose $k \ge 7$ and the result is true for all values less than $k$, and let $G$ be a non-planar cyclically $4$-connected cubic graph with $2k$ vertices other than $V_{2k}$. By Theorem \ref{Wormald}, $G$ is obtained from a cyclically $4$-connected cubic graph $G'$ with $2(k-1)$ vertices by bridging a pair of edges with cycle spread at least $(1,1)$. By Lemma \ref{exchange-corollary}, one of the following is true: 

\begin{enumerate}
\item $G' \cong Q_{2(k-1)}$, $G \cong Q_{2k}$, and $G$ is obtained from $G'$ by bridging a pair of edges with cycle spread $(1,1)$; a contradiction since $G$ is non-planar.

\item $G' \cong V_{2(k-1)}$, $G \cong V_{2k}$ and $G$ is obtained from $G'$ by bridging a pair of edges with cycle spread $(1,1)$; a contradiction since $G \not \cong V_{2k}$.

\item The bridge that produces $G$ from $G'$ can be exchanged for a bridge of a pair of edges with cycle spread at least $(1,2)$.

\end{enumerate}

In other words, the third case of Lemma \ref{exchange-corollary} holds. If the graph $G' \not \cong V_{2(k-1)}$, we are done. Otherwise, suppose $G' \cong V_{2(k-1)}$. Since $k \ge 7$ and $G \not \cong V_{2k}$, there must be another edge in $G$ that can be unbridged to obtain a graph $G'' \not \cong V_{2(k-1)}$ with $2(k-1)$ vertices. By the induction hypothesis, $G''$ can be obtained from $Q_8$ by bridging pairs of edges with cycle spread $(1,2)$. By Lemma \ref{exchange-corollary}, since $G \not \cong V_{2k}$ and $G'' \not \cong V_{2(k-1)}$, $G$ is obtained from $G''$ by bridging a pair of edges with cycle spread at least $(1,2)$. The result follows.
\end{proof}

Our goal is to generate precisely the objects we are looking for, in this case non-planar cyclically $4$-connected cubic graphs, and Theorem \ref{MainTheorem-2} accomplishes this.

The next result is similar to Theorem \ref{MainTheorem-2} but it for planar graphs and the process begins with some $Q_{2k}$, not simply $Q_8$.

\begin{theorem} \label{MainTheorem-3} Let $G$ be a planar cyclically $4$-connected cubic graph with $n \ge 8$ vertices, where $n$ is even, and $G \not \cong Q_n$. Then $G$ is obtained from $Q_{2k}$ for some $k \ge 4$ by bridging pairs of edges with cycle spread at least $(1, 2)$.
\end{theorem}

\begin{proof} 

The proof is by induction on $n$. Since $Q_8$ and $Q_{10}$ are the only planar cyclically $4$-connected cubic graphs of size $8$ and $10$, respectively, the result is true for graphs of those sizes. The only planar cyclically $4$-connected cubic graph of size $12$ other than $Q_{12}$ is the graph shown in the lower right-hand portion of Figure \ref{exchange-figure}, which is obtained from $Q_{10}$ by bridging a pair of edges with cycle spread $(1,2)$. Thus the result holds for $n= 8, 10, 12$.

Suppose $n \ge 14$ and the result is true for all values less than $n$, and let $G$ be a planar cyclically $4$-connected cubic graph with $n$ vertices other than $Q_n$. By Theorem \ref{Wormald}, $G$ is obtained from a smaller cyclically $4$-connected cubic graph $G'$ by bridging a pair of edges with cycle spread at least $(1,1)$. In other words, $G'$ is obtained from $G$ by an unbridging operation and since the unbridging operation preserves planarity, $G'$ must be planar too. By Lemma \ref{exchange-corollary}, one of the following is true:

\begin{enumerate}

\item $G' \cong Q_{n-2}$ and $G \cong Q_n$; a contradiction to the hypothesis of the theorem.

\item $G' \cong V_{n-2}$ and $G \cong V_n$; a contradiction since $G$ is planar.

\item $G$ is obtained from $G'$ by bridging a pair of edges with cycle spread at least $(1,2)$.

\end{enumerate}

By the induction hypothesis, $G'$ is obtained from $Q_{2k}$ for some $k \ge 4$ by bridging pairs of edges with cycle spread at least $(1,2)$. Since $G$ is obtained from $G'$ by bridging a pair of edges with cycle spread at least $(1,2)$, the result follows.
\end{proof}


\section{Applications}

In order to construct cyclically $4$-connected cubic graphs using Theorem \ref{Wormald} we must identify pairs of edges with cycle spread at least $(1,1)$. In order to construct cyclically $4$-connected cubic graphs using Theorems \ref{MainTheorem-2} and \ref{MainTheorem-3}, we must identify pairs of edges with cycle spread at least $(1,2)$. These procedures are exhaustive. Moreover, if we can identify pairs of edges with cycle spread at least $(2,2)$, then beginning with cyclically $5$-connected cubic graphs we can bridge these pairs to form larger cyclically $5$-connected cubic graphs, although this procedure is not exhaustive.

Our approach to identifying a pair of edges $ab$ and $cd$, as shown in Figure \ref{BG-2}, with a specified minimum cycle spread is as follows:

\begin{enumerate}

\item Choose a pair of non-adjacent vertices $b$ and $c$ in a cyclically $t$-connected cubic graph $G$, where $t$ is $4$ or $5$; and

\item Identify edges by strategically selecting neighbors of each of them. 

\begin{itemize}

\item[a)] If we choose a vertex $d \in N(c) \backslash \{b\}$ and a vertex $a \in N(b) \backslash \{c, d\}$ then the edges $ab$ and $cd$ will have cycle spread at least $(1,1)$ because $ab$ and $cd$ will be non-adjacent. In other words, if vertices $b$ and $c$ are non-adjacent, and we choose vertex $d$ to be a neighbor of $c$, and choose vertex $a$ to be a neighbor of $b$ avoiding $d$, then edges $ab$ and $cd$ will have cycle spread at least $(1,1)$.

\item[b)] If we choose a vertex $d \in N(c) \backslash (N(b) \cup \{b\})$ and a vertex $a \in N(b) \backslash \{c, d\}$, then $ab$ and $cd$ will have cycle spread at least $(1,2)$ because $a \not \in \{b, c\}$ and $b \not \in N(c) \cup N(d)$. In other words, if vertices $b$ and $c$ are non-adjacent, and we choose vertex $d$ to be a neighbor of $c$ avoiding neighbors of $b$, and we choose vertex $a$ to be a neighbor of $b$ avoiding $d$, then edges $ab$ and $cd$ will have cycle spread at least $(1,2)$.

\item[c)] If we choose a vertex $d \in N(c) \backslash (N(b) \cup \{b\})$ and a vertex $a \in N(b) \backslash (\{c, d\} \cup N(c) \cup N(d))$, then $ab$ and $cd$ will have cycle spread at least $(2,2)$ because $a, b \not \in \{c, d\} \cup N(c) \cup N(d)$. In other words, if vertices $b$ and $c$ are non-adjacent, and we choose vertex $d$ to be a neighbor of $c$ avoiding neighbors of $b$, and we choose vertex $a$ to be a neighbor of $b$ avoiding neighbors of $c$, then $ab$ and $cd$ will have cycle spread at least $(2,2)$.

\end{itemize}
\end{enumerate}

\begin{figure}[h]
\centering
\includegraphics[width=2.3in]{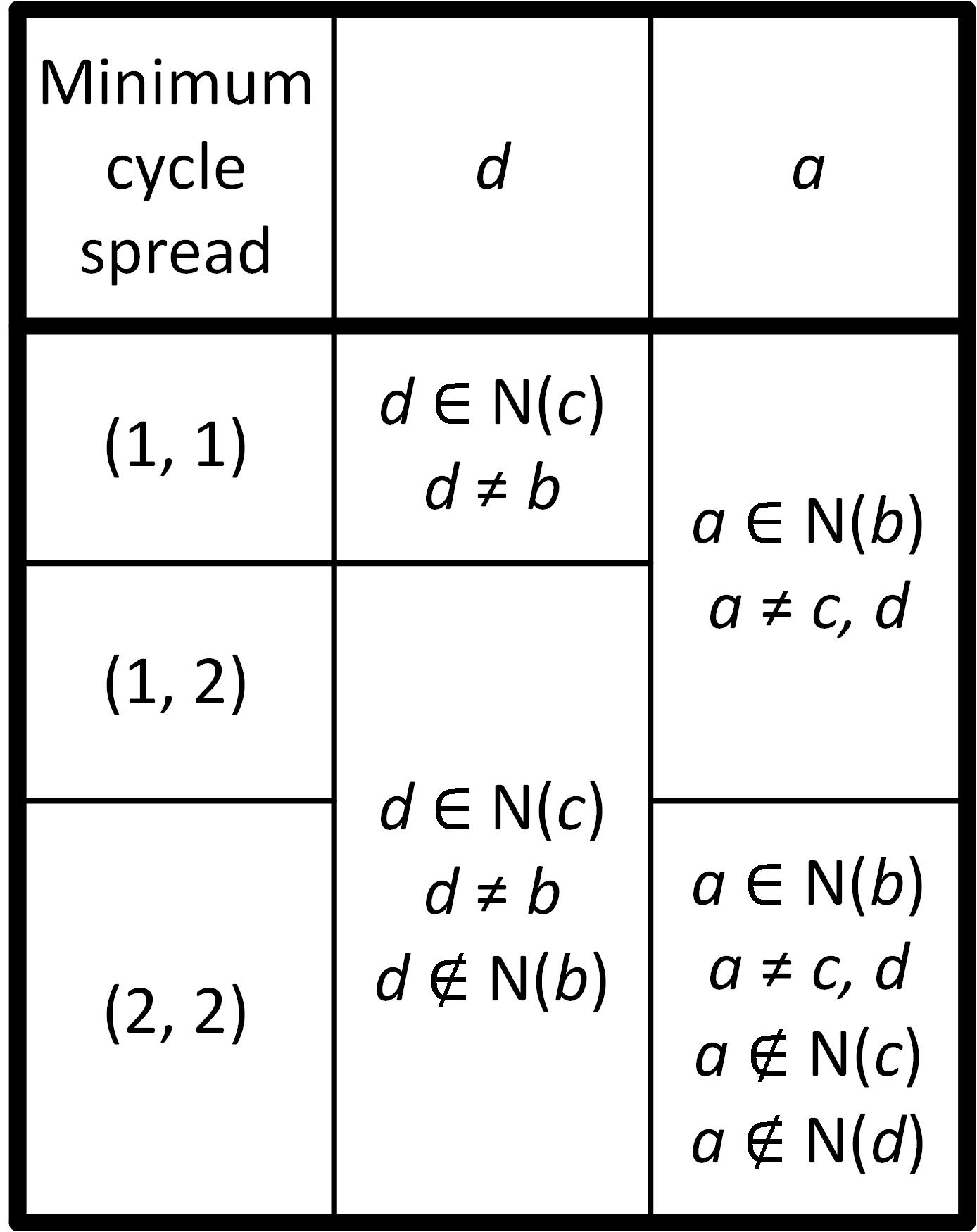}
\caption{Criteria for selection of vertices $d$ and $a$ for various minimum cycle spreads. \label{cycle-spread-criteria-table}}
\end{figure}

The criteria for selecting $d$ and $a$ are summarized in the table in Figure \ref{cycle-spread-criteria-table}. As the table indicates, from cycle spread $(1,1)$ to cycle spread $(1,2)$ the criteria for selecting $d$ become more restrictive but the criteria for selecting $a$ do not. From cycle spread $(1,2)$ to cycle spread $(2,2)$ the criteria for selecting $d$ do not change, but the criteria for selecting $a$ become more restrictive.

Using this strategy, we can identify pairs of edges to bridge in order to implement Wormald's results for propagating cyclically $4$- or $5$-connected graphs, or our results for propagating cyclically $4$-connected graphs using a smaller set of bridge additions. We illustrate the procedure by comparing  bridge additions obtained of pairs of edges with cycle spread at least $(1,1)$ or $(1,2)$ from a pair of non-adjacent vertices in $Q_8$. 

\begin{figure}[h]
\centering
\includegraphics[width=5.15in]{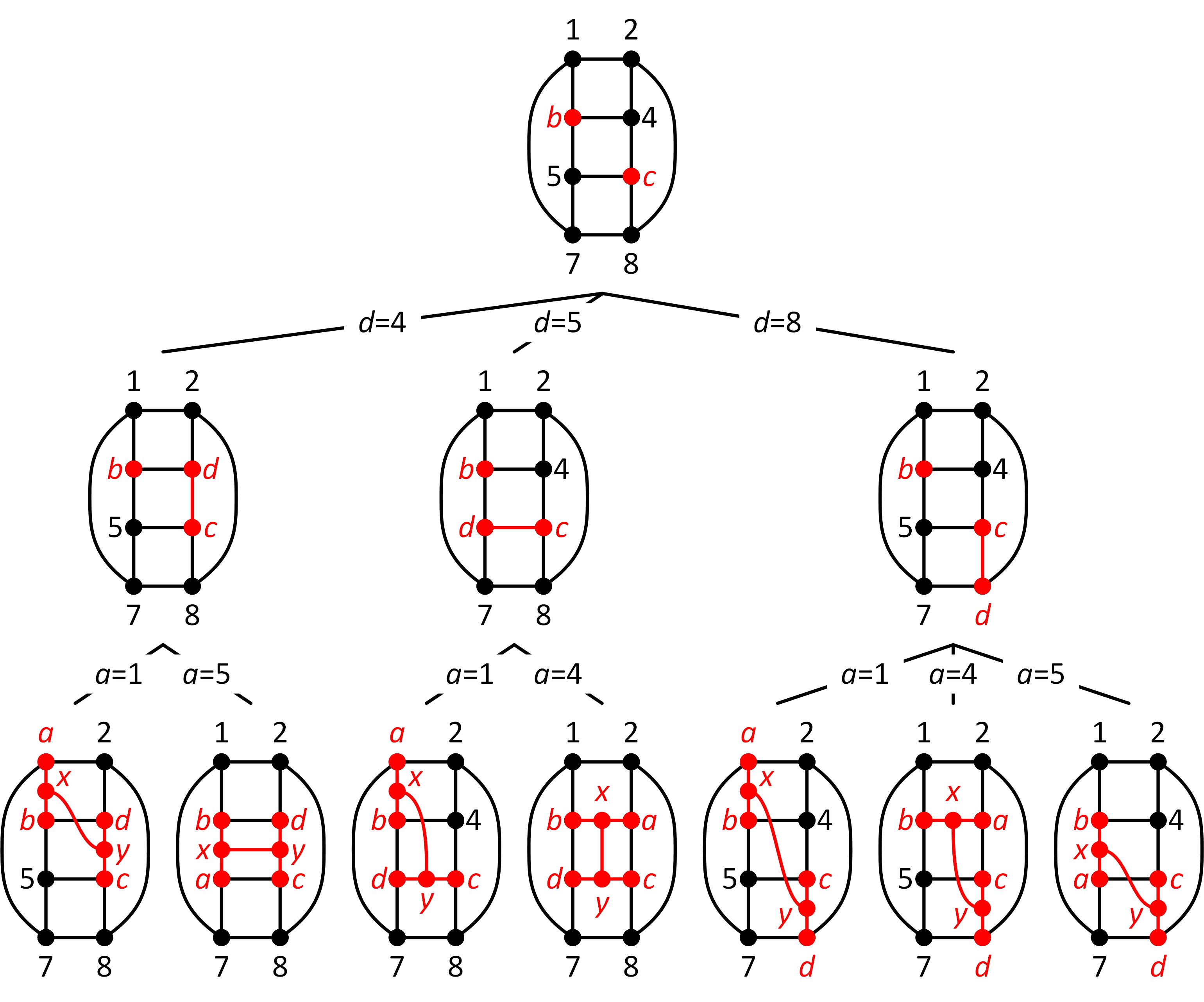}
\caption{Cyclically $4$-connected cubic graphs obtained by bridging pairs of edges with cycle spread at least $(1,1)$ in $Q_8$, using the pair of non-adjacent vertices $b$ and $c$. The graphs under $d=8$ are the ones that are checked when seeking pairs of edges with cycle spread at least $(1,2)$. \label{q8-1-1-all-figure}}
\end{figure}

The graph at the top of Figure \ref{q8-1-1-all-figure} is $Q_8$ with two non-adjacent vertices $b$ and $c$ indicated. Following (a) in the procedure above for cycle spread at least $(1, 1)$, there are three choices for $d$, namely vertices 4, 5, and 8, the neighbors of $c$. The middle row of graphs in the figure shows $Q_8$ with edge $cd$ indicated for each choice of $d$.

Next, we need to select a neighbor $a$ of $b$ that is not $d$.  

\begin{itemize}

\item[i)] If $d = 4$ then the choices are $a=1$ or $a=5$. The two graphs on the left-hand side of the bottom row of Figure \ref{q8-1-1-all-figure} show the graphs resulting from each of these choices. The graph obtained with $d=4$ and $a=1$ is $A_{10,1}$ and the graph obtained with $d=4$ and $a=5$ is $Q_{10}$.

\item[ii)] If $d = 5$ then the choices are $a=1$ or $a=4$. The resulting bridge additions are shown in the center of the bottom row of Figure \ref{q8-1-1-all-figure}. The graph obtained with $d=5$ and $a=1$ is $A_{10,1}$ and the graph obtained with $d=5$ and $a=4$ is $Q_{10}$.

\item[iii)] If $d=8$ then any neighbor of $b$ may be selected as $a$, so $a=1$, $a=4$, or $a=5$. The resulting bridge additions are shown to the right of the bottom row of Figure \ref{q8-1-1-all-figure}. The graph obtained with $d=8$ and $a=1$ is $A_{10,2}$. The graphs obtained with $d=8$ and either $a=4$ or $a=5$ are both $A_{10,1}$.

\end{itemize}

Next we will apply method (b) in the procedure above to identify pairs of edges with cycle spread at least $(1,2)$ from the vertices $b$ and $c$ in $Q_8$. Observe that we can obtain the same non-planar results while generating fewer graphs. Since vertices 4 and 5 are neighbors of $b$, the only neighbor of $c$ that is eligible for selection as $d$ is vertex 8, which is the same as case (iii) above. Once we select $d=8$, since it is not a neighbor of $b$ itself, any neighbor of $b$ may be selected as $a$, namely vertex 1, 4, or 5. The three graphs obtained are shown to the right in the bottom row of Figure \ref{q8-1-1-all-figure}. The choices of $a$ with $d=8$ are exactly the same for cycle spread at least $(1,1)$ and $(1,2)$. We obtain the cyclically $4$-connected non-planar cubic graphs $A_{10,1}$ and $A_{10,2}$, but not the cyclically $4$-connected planar graph $Q_{10}$.

\begin{figure}[h]
\centering
\includegraphics[width=2.65in]{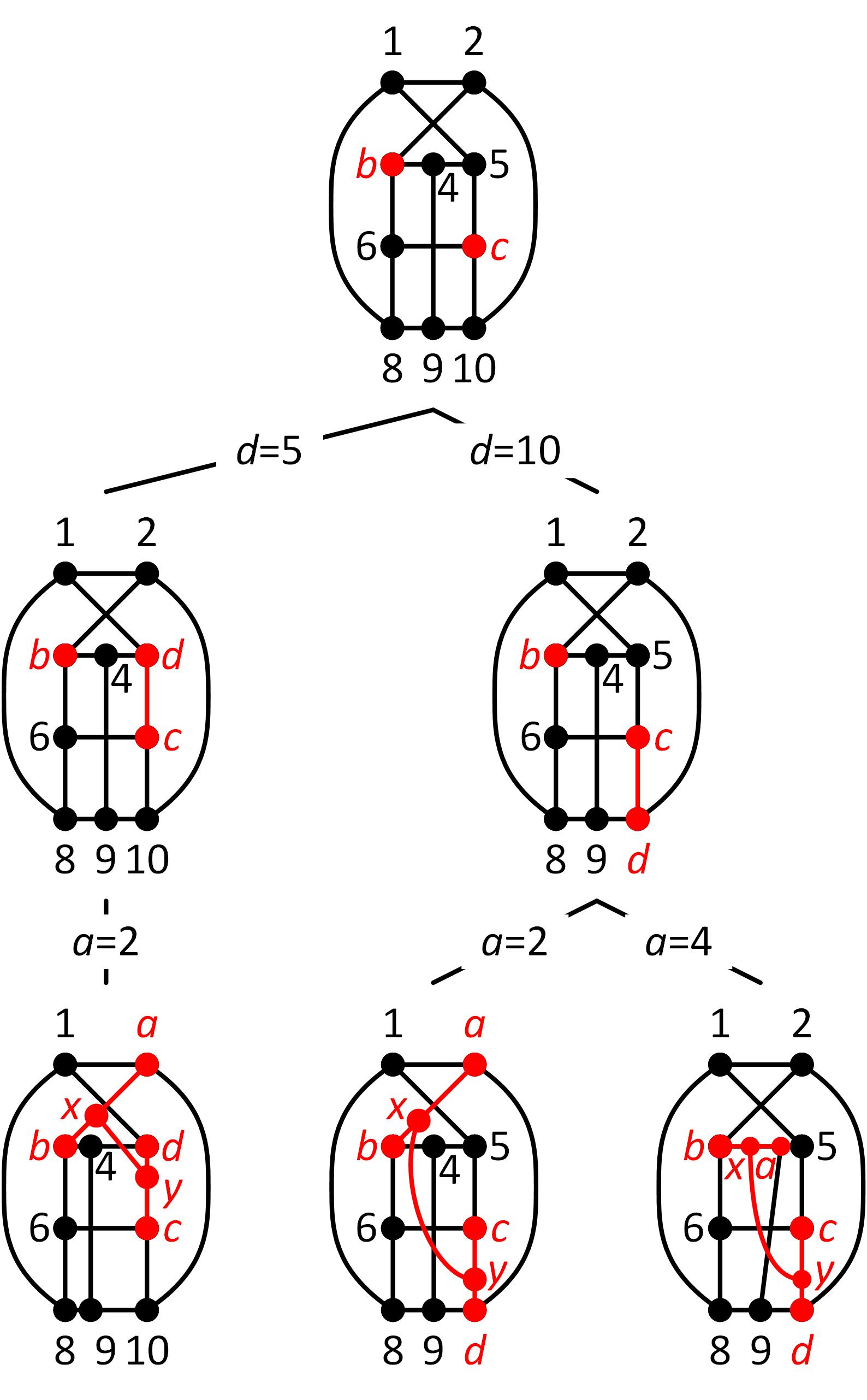}
\caption{Cyclically $5$-connected cubic graphs obtained by bridging pairs of edges with cycle spread at least $(2, 2)$ in $P_{10}$, using the pair of non-adjacent vertices $b$ and $c$. \label{p10-2-2-figure}}
\end{figure}

Finally, Figure \ref{p10-2-2-figure} shows the bridge additions obtained of pairs of edges with cycle spread at least $(2,2)$ from a pair of non-adjacent vertices in $P_{10}$. Specifically, we use method (c) in the procedure to identify pairs of edges with cycle spread at least $(2,2)$ from vertices $b$ and $c$. The candidates for $d$ are the neighbors of $c$ that are not also $b$ or neighbors of $b$, namely $d=5$ and $d=10$. The middle row of the figure shows $P_{10}$ with each choice for edge $cd$ indicated. Vertex $a$ must be a neighbor of $b$ that is not $d$, $c$, or any neighbor of $c$ or $d$. If $d=5$ the only choice is $a=2$; the resulting graph is shown on the left in the bottom row of Figure \ref{p10-2-2-figure}. If $d=10$ then there are two choices, $a=2$ or $a=4$. The resulting graphs are shown to the right in the bottom row of Figure \ref{p10-2-2-figure}.

We developed an algorithm based on these theorems and an implementation in Python for generating non-planar cyclically $4$-connected cubic graphs, by translating bridge additions into edge additions and vertex splits. This was an extension to the implementation of an algorithm to generate minimally $3$-connected graphs described in \cite{Kingan2021}, which uses the certificate generating method in McKay's nauty system \cite{McKay1981} to eliminate isomorphic duplicates from its output. Early checking for isomorphic duplicates at each stage of the implementation considerably reduces the overall number of certificate generations that must be done.

Using this program, all non-planar cyclically 4-connected cubic graphs up to 20 vertices were generated on a PC with an Intel Core I5-4460 CPU at 3.2 GHz and 16 Gb of RAM, in around 35 minutes. The graphs were validated using a separate procedure which directly checked for size 3 edge cuts with cycles in the remaining components, and the number of graphs generated was compared to the published counts in the Online Encyclopedia of Integer Sequences.



\end{document}